\documentclass[a4paper,12pt]{article}

\usepackage{amsfonts}
\usepackage{amsthm}
\usepackage{amsmath}
\usepackage{amsfonts}
\usepackage{latexsym}
\usepackage{amssymb}
\newtheorem{fed}{Definition}[section]
\newtheorem{teo}[fed]{Theorem}
\newtheorem{lem}[fed]{Lemma}
\newtheorem{cor}[fed]{Corollary}
\newtheorem{pro}[fed]{Proposition}
\theoremstyle{definition}
\newtheorem{rem}[fed]{Remark}
\newtheorem{rems}[fed]{Remarks}
\newtheorem{exa}[fed]{Example}

\def\inv{^{-1}}

\def\beq{\begin{equation}}
\def\eeq{\end{equation}}
\def\N{\mathbb{N}}
\def\F{\mathcal{F}}
\def\G{\mathcal{G}}

\def\cB{\mathcal{B}}

\def\cH{\mathcal{H}}
\def\cK{\mathcal{K}}

\def\ese{\mathcal{M}}
\def\ete{\mathcal{T}}
\def\eme{\mathcal{M}}
\def\ene{\mathcal{N}}

\def\cW{\mathcal{W}}

\def\inc{\subseteq}
\def\sudi{\, \dot {+} \, }

\def\bdem{\begin{proof}}
\def\edem{\end{proof}}

\def\orto{^\perp}
\def\inc{\subseteq}
\def\sii{ if and only if }
\def\inv{^{-1}}
\def\*A{\#\sb A}

\def\eps{\varepsilon}

\def\H{{\cal H}}

\def\glh{Gl(\cH)}
\def\ca{L(\H ) }

\def\cH{{\cal H}}

\def\cM{{\cal M}}

\def\rai{^{1/2}}
\def\mrai{^{-1/2}}

\def\api{\langle}
\def\cpi{\rangle}

\def\noi{\noindent}

\def\bm{\left(\begin{array}}
\def\em{\end{array}\right)}
\def\ben{\begin{enumerate}}
\def\een{\end{enumerate}}
\def\barr{\begin{array}}
\def\earr{\end{array}}
\def\iiff{if and only if }
\def\inv{^{-1}}

\def\H{{\cal H}}
\def\lh{{L(\H)}}
\def\lh+{{\lh^+}}

\def\la{\lambda}
\def\eps{\varepsilon}

\def\EOE{\hfill$\triangle$}

\DeclareMathOperator*{\convsotipre}{\nearrow}

\newcommand{\pint}[1]{\displaystyle \left \langle #1 \right\rangle}

\newcommand{\hil}{\mathcal{H}}
\newcommand{\op}{L(\mathcal{H})}

\newcommand{\Dpint}[1]{\displaystyle \left \langle #1 \right\rangle_D}
\newcommand{\gen}[1]{\mbox{span}\left\{#1\right\}}
\newcommand{\cl}[1]{\mbox{cl}\left(#1\right)}
\newcommand{\gencl}[1]{\overline {\mbox{span}}\left\{#1\right\}}

\newcommand{\convk}{\xrightarrow[k\rightarrow\infty]{}}

\newcommand{\convsot}{\xrightarrow[n\rightarrow\infty]{\mbox{\tiny{SOT}}}}

\newcommand{\convsoti}{\convsotipre_{n\rightarrow\infty}^{\mbox{\tiny
\textbf SOT}}}

\newcommand{\bol}[1]{#1_1}
\newcommand{\fram}{\{f_n\}_{n\in \mathbb{N}}}

\newcommand{\framg}{\{g_n\}_{n\in \mathbb{N}}}

\newcommand{\bon}{\{e_k\}_{k\in \mathbb{N}}}
\newcommand{\bonn}{\{e_n\}_{n\in \mathbb{N}}}

\newcommand{\deledos}{\mathcal{D}(\ell^2)}

\newcommand{\pas}[2]{P_{#1,\ #2}}

\def\Chi{\mbox{\Large{$\chi$}}}
\def\noi{\noindent}
\def\QED{\hfill $\square$}

\newcommand{\dist}[2]{\mbox{d}\left(#1,\, #2\right)}

\newcommand{\angf}[2]{c\left[\,#1,\,#2\,\right]}
\newcommand{\sang}[2]{s\left[\,#1,\,#2\,\right]}
\begin{document}
% ----------------------------------------------------------------
%\date{}

% Title ----------------------------------------------------------

\title{{\Large\textbf{Nullspaces and  frames}
\thanks{Partially supported by CONICET (PIP 2083/00), UBACYT I030, UNLP (11 X350) 
and ANPCYT (PICT03-9521)}
}
%\footnote{2000 AMS-MSC: Primary 46C07, 47A62, 46C05. 
%} \footnote{
%{\bf Keywords:} shorted matrix, spectral order, positive matrices.
%}
}
\author{J. Antezana
%\thanks{Partially supported by UNLP (11 X350)} 
, G. Corach
%\thanks{Partially supported by CONICET (PIP 2083/00), UBACYT I030 and ANPCYT (PICT03-9521)} , 
, M. Ruiz 
%\thanks{Partially supported by UNLP (11 X350) and ANPCYT (PICT03-9521)}
and D. Stojanoff
%\thanks{Partially supported by CONICET (PIP 2083/00), UBACYT I030, UNLP (11 X350) 
%and ANPCYT (PICT03-9521)}
}
%\author{Jorge Antezana, Gustavo Corach, Mariano Ruiz and Demetrio
%Stojanoff}
%\date{August 15, 2003}

%\thanks
%\keywords {Oblique projections, operator ranges, positive operators}

\maketitle

\vglue.3truecm

\noi
{\bf Jorge Antezana, Mariano Ruiz and Demetrio Stojanoff }

\noi
Depto. de Matem\'atica, FCE-UNLP,  La Plata, Argentina
and IAM-CONICET

\noi
e-mail: antezana@mate.unlp.edu.ar, mruiz@mate.unlp.edu.ar\\ and demetrio@mate.unlp.edu.ar

\medskip

\noi
{\bf Gustavo Corach (corresponding author)}

\noi
Depto. de Matem\'atica, FI-UBA  and IAM-CONICET, \\
Saavedra 15,  Piso 3  (1083),\\
Ciudad Aut\'onoma de Buenos Aires, 
Argentina.\\
Phone: (54) (11) 4954 - 6781 \\
Fax: (54) (11) 4954 - 6782 \\
e-mail: gcorach@fi.uba.ar

\medskip
\noi

\medskip
\noi

\vglue1truecm
\begin{abstract}
In this paper we give new characterizations of Riesz and 
conditional Riesz frames in terms of the properties of the nullspace 
of their synthesis operators. On the other hand, we also 
study the oblique dual frames whose coefficients in the 
reconstruction formula minimize different weighted norms.
\end{abstract}

%\newpage

%\tableofcontents

\section {Introduction}

Frames were introduced by Duffin and Schaeffer \cite{[DuSh]} in the 
context of nonharmonic Fourier series, and they have been intensively applied in 
wavelet and frequence analysis theories since the work of Daubechies, 
Grossmann and Meyer \cite{[DGM]}. 
%They play an important role in many applications .
Today, frame-like expansions are fundamental in a wide range of 
disciplines (see for example \cite{[DuSh]}, \cite{[HeiWal]} or \cite{[Y]}), including the analysis and 
design of oversampled filter banks and error corrections codes. 
%Wavelet families have been used in quantum mechanics and many other areas of theoretical physics.

A frame is a redundant set of vectors in a Hilbert space that leads to expansions of 
vectors (signals) in terms of the frame elements. More precisely, a sequence of 
vectors $\F = \fram$ in a (separable) Hilbert space $\H$ is a \textsl{frame} (for $\hil$) if
there exist numbers $ A,B>0$ such that, for every $f \in \H$,
\beq\label{frame1}
A\|f\|^2\leq \sum_{n\in \mathbb{N}}
|\pint{f , f_n}|^2 \leq B\|f\|^2 \ .
\end{equation}
Associated with each frame there exists an operator $T:\ell^2\to\hil$ 
defined by $T(e_n)=f_n$, where $\cB=\bonn$ denotes the canonical basis 
of $\ell^2$; $T$ is called the synthesis operator of $\F$.

The results of this paper can be divided in two parts. The main 
results of the first part are devoted to the 
study of Riesz frames and conditional Riesz frames through the structure and 
geometric properties of the nullspace of their synthesis 
operators. Riesz and conditional Riesz frames were introduced 
by Christensen in \cite{[Chri3]} (see definitions in section 3). 
These frames are important because they behave well with respect to the projection method.  
In general, frame theory describes how to choose the corresponding coefficients  to expand a given 
vector in terms of the frame vectors. However, in applications, to obtain these 
coefficient requires the inversion of an operator on $\hil$. The projection method was 
introduced by Christensen in \cite{[Chri1]}  to avoid this problem. We refer the 
interested reader to \cite{[CazChri3]}, \cite{[Chri1]}, \cite{[Chri3]} or \cite{[liChr]} 
for more information about the projection method.
In \cite{[ACRS]} we found a characterization of Riesz frames by
studying the nullspace of the synthesis operator.
Namely, if  the nullspace $N(T)$ has a certain geometric property of
compatibility with the  closed subspaces spanned by subsets of $\cB$,
then $\F $ is a Riesz frame, and conversely.
In section 3, we extend these results for conditional Riesz frames 
and give some new characterizations in terms of angles.

Throughout the second part of this work we study the so-called oblique dual frames. 
Let $\fram$ be a frame for the closed subspace $\cW\subseteq \hil$, and 
let $\ese\subseteq \hil$ be another closed subspace such that 
$\hil=\cW\sudi \ese^\bot$ ($\sudi$ means a non necessarily orthogonal direct sum). 
The sequence $\framg$ in $\ese$ is an oblique dual 
frame of $\fram$ (see Li \cite{[Li]} or Li and Ogawa \cite{[LO1]} and \cite{[LO2]}) if
\[
f=\sum_{n=1}^\infty\pint{f,\ g_n}f_n \hspace{1cm} \forall \ f\in\cW.
\]
Among the oblique dual frames, there exists a particular class with the 
minimal norm property. Recall that a dual frame $\framg$ has the minimal 
norm property if the coefficients $\{\pint{f,\ g_n}\}_{n\in\N}$ that 
appear in the reconstruction formula have minimal $\ell^2$ norm. 
	
If $\cB=\bonn$ denotes the canonical orthonormal basis for $\ell^2$ and $T$ 
is the synthesis operator of $\fram$, then Christensen and Eldar \cite{[ChriEl]} proved that the 
minimal norm oblique dual frames have the form
\[
\framg=\left\{B(T^*B)^\dagger e_n\right\}_{n\in\N} \ ,
\]
where $B$ is any bounded operator with $R(B)=\ese$. From 
the point of view of sampling theory, the operator $B$ can be interpreted 
as the synthesis operator associated to the frame used to sample the signals.

In this work, we are interested in duals frames which 
lead to reconstruction coefficients that have minimal norm, but with 
respect to some weighted norms. Recall that weighted norms 
in $\ell^2$ arise from inner products obtained by perturbing the 
original one with invertible positive operators which are diagonal in the canonical basis.
In section 4 we give explicit formulae for dual frames which 
minimize a given weighted norm, and we prove that in the 
case of Riesz frames, if the sampling frame is fixed, 
then the norms of the preframes operators corresponding to the dual frames 
which minimize the different weighted norms are uniformly bounded from above.

We thank Ole Christensen for his useful comments.

\section{Preliminaries}

Let $\hil$ be a separable Hilbert space and $\ca$ the
algebra of bounded linear operators on $\hil$. $\glh$ 
denotes the group of invertible operators in $\ca$, and
$\glh^+$ the set of positive definite invertible operators on $\H$.
For an operator $A
\in \ca$, $R(A)$ denotes the range of $A$, $N( A)$ the
nullspace of $A$, $\sigma (A)$ the spectrum of $A$, $A^*$ the adjoint
of $A$, $\rho(A)$ the spectral radius of $A$ and $\|A\|$ the
operator norm of $A$; if $R(A)$ is closed, $A^\dagger$ is the
Moore-Penrose pseudoinverse of $A$. We use the fact  that $A$ is an isometry 
(resp. coisometry) if $A^*A = I$ (resp. $AA^*=I$).
Given a closed subspace $ \ese $ of $ \hil$, $P_\ese$ denotes the
orthogonal (i.e., selfadjoint) projection onto $\ese$. If $B \in
\ca $ satisfies $P_\ese B P_\ese = B$, we consider the compression
of $B$ to $\ese$, (i.e., the restriction of $B$ to $\ese$, which  is an
operator on $\ese$), and we say that we
consider $B$ as $acting$ on $\ese$. 
Given a subspace $\ese$ of $\hil$, its unit ball is denoted by
$\bol{\ese}$, and its closure by $\overline{\ese}$ or $\cl{\ese}$.
If $\ene$ is another subspace of $\hil$, we denote $\ese \ominus \ene :=
\ese \cap \ene \orto$. If $\ese \cap \ene = \{0\} $, we denote by $\ese \sudi \ene$ the 
(direct) sum of the two subspaces. If the sum is orthogonal, we write $\ese \oplus \ene$.
 The distance between two subsets $\eme$ and $\ene$ of $\hil$ is
$
\dist{\eme}{\ene}=\inf\{\|x  -y  \|:\; x  \in \eme \;\; y  \in \ene\}.
$

\subsection{Angle between closed subspaces}
We shall recall the definition of angle between closed subspaces of
$\hil$. We refer the reader to the nice survey  of Deutsch
\cite{[De]} and the books by Kato \cite{[Kato]} and  Havin and
J\"oricke \cite{[HJ]} for details and proofs.

\begin{fed}\rm
Given two closed subspaces $\eme$ and $\ene$, let
$\tilde \ene = \ene \ominus (\eme \cap \ene ) $ and
$\tilde \eme = \eme \ominus (\eme \cap \ene ) $.
The \textbf{angle} between $\eme$ and
$\ene$ is the angle in $[0,\pi/2]$ whose cosine
is 
\[
\angf{\eme}{\ene}=\sup\{\,|\pint{x,\,y}|:\;x\in \tilde \eme ,
\;y\in \tilde \ene \;\mbox{and}\;\|x\|=\|y\|=1 \}
\]
The $sine$ of this angle is denoted by $\sang{\eme }{\ene }$.
\end{fed}

\bigskip

\noi Now, we state some known results concerning angles and closed
range
 operators (see  \cite{[De]}).

\begin{pro}\label{propiedades elementales de los angulos}
Let $\eme$ and $\ene$ be two closed subspaces of $\hil$.
Then \ben
\item $\angf{\eme}{\ene}=\angf{\ene}{\eme}=\angf{\tilde \eme}{\ene}=
  \angf{\eme}{\tilde \ene}$.
\item $\angf{\eme}{\ene}<1$   \sii $\eme+\ene$ is closed.
\item $\angf{\eme}{\ene}= \angf{\eme^\bot}{\ene^\bot}$
\item $\angf{\eme}{\ene}=\|P_\eme P_{\tilde \ene}\|=\|P_{\tilde \eme}P_\ene \|
=\|P_\eme P_\ene P_{(\eme\cap \ene)^\bot}\|=\|P_\eme P_\ene-P_{\eme\cap \ene}\|$
\een
\end{pro}

\smallskip
\begin{pro}[Bouldin \cite{[Bou]}; see also \cite{[De]}]\label{producto con rango cerrado}
Let $A,B\in\op$ such that $R(A)$ and $R(B)$ are closed. Then, $AB$ has closed range \sii \
$\angf{R(B)}{N(A)}<1$.
\end{pro}

\smallskip
\begin{pro}[Kayalar-Weinert \cite{[Kayalar]}; see also \cite{[De]}]
\label{aproximacion a la interseccion}
Let $P$ and $Q$ two orthogonal projections defined on $\hil$.
Then,
\[
\|(P Q)^k-P\wedge Q\|=\angf{R(P)}{R(Q)}^{2k-1}
\]
where $P\wedge Q$ is the orthogonal projection onto $R(P)\cap R(Q)$.
\end{pro}

\noi
Finally,  we give a characterization of $\sang{\eme}{\ene}$ in terms
of distances:

\begin{pro}\label{seno}
Let $\eme$ and $\ene$ be to closed subspaces of $\hil$.
Denote $\tilde \ene = \ene \ominus (\eme \cap \ene ) $ and
$\tilde \eme = \eme \ominus (\eme \cap \ene ) $. Then
$$
\sang {\eme}{\ene}= \dist {\tilde\eme_1 }{\ene} =\dist {\tilde \ene_1 }{\eme}.
$$
\end{pro}
\bdem
By Proposition \ref{propiedades elementales de los angulos},
we can suppose that $\eme\cap \ene =\{0\}$,
i.e., $\eme = \tilde \eme$. By the definition of the sine and
Proposition \ref{propiedades elementales de los angulos},
$\sang{\eme}{\ene}^2= 1-\|P_\eme P_{\ene}\|^2 $. On the other hand, as $\dist{x}{\ene} = \|P_{\ene\orto } \, x\|$ for every $x \in \H$, we have that
$$
\barr{rl}
\dist {\eme_1 }{\ene}^2 &= \inf \{ \|P_{\ene\orto }\,x\|^2 : x\in \eme_1  \}%\\&\\
= \inf \{1-  \|P_{\ene }\,x\|^2 : x\in \eme_1  \}\\&\\
&= 1-\sup \{ \|P_{\ene }\,x\|^2 : x\in \eme_1  \} = 1-\|P_\ene P_{\eme}\|^2=
1-\|P_\eme P_{\ene}\|^2.
\earr
$$
\edem

\subsection{The reduced minimum modulus}
\noi
\begin{fed}\rm
The \textit{reduced minimum modulus} $\gamma(T)$ of an operator
$T\in L(\hil)$ is defined by
\begin{equation}\label{gamma}
\gamma (T)= \inf \{ \|Tx\| :  \|x\|=1 \; , \;  x\in N(T)^\bot \}
\end{equation}
It is well known that $\gamma (T)=\gamma (T^*) =\gamma (T^*T)\rai $. Also, it can be
shown that an operator $T$ has closed range if and only if $\gamma
(T)>0$.  In this case, $  \gamma (T)= \|T^\dag \|^{-1} $. 
\end{fed}

\medskip

\noi
The following result is an easy consequence of equation (\ref{gamma}):

\begin{lem}\label{L:inversible}
Let $B\in \ca $ with $B$ invertible. Then,
\[
\|B^{-1}\|^{-1} \gamma (T) \leq \gamma (BT) \leq \|B\| \gamma
(T)  .
\]
Moreover, the same formula follows, replacing $\|B^{-1}\|^{-1} $ by $\gamma(B)$,
if  $R(B)$ is closed and $R(T) \inc N(B)\orto$.
\end{lem}

\begin{lem}\label{L:coisometria}
Let $T \in \ca$ be a  partial isometry (i.e., $TT^*$ is a projection),
$\eme$ a closed subspace of
$\hil$ and $P_\eme$ the orthogonal projection onto $\eme$.  Then
$$
\gamma (TP_\eme)=\sang{N(T)}{ \eme} \ .
$$
\end{lem}

\begin{proof}
\def\ere{\mathcal{R}}
Denote $\ene = N( T)$ and $\ere = \ene\orto$.
Since $T$ acts isometrically on $\ere$,  it is clear by equation (\ref{gamma}) that
$$\gamma (TP_\eme) = \gamma (TP_\ere P_\eme)=\gamma (P_\ere P_\eme) .
$$
Since $N(P_\ere P_\eme ) = \eme^\bot \oplus (\eme\cap \ene) $, it follows that
$N(P_\ere P_\eme)^\bot  = \eme\cap (\eme\cap \ene)^\bot = \tilde \eme$.
Then, by Proposition \ref{seno},
$$
\gamma (P_\ere P_\eme) = \inf _{x\in \tilde \eme_1} \|P_\ere x\| =
\inf _{x\in \tilde \eme_1 } \dist{x}{ \ene}   =
\dist {\tilde \eme_1 }{\ene} = \sang{\ene}{ \eme} \ .
$$
\end{proof}

\noi
The next result was proved in \cite {[ACRS]}. We include
a short proof for the sake of completeness.

\begin{pro}\label{P:gammas}%(cita a la proposicion en weighted)
If $T \in \ca $ has closed range and  $\eme$ is a closed
subspace of $\hil$ such that $ \angf{N(T)}{ \eme}<1$ (so that $TP_\eme$ has closed
range), then
\begin{equation}\label{desigualdad}
 \gamma (T)\ \sang{N(T)}{ \eme} \leq \gamma (TP_\eme)\leq \|T\| \ \sang{N(T)}{ \eme} .
\end{equation}
\end{pro}

\begin{proof}
Take $B=|T^*|=(TT^*)\rai$. It is well known that $R(B) = R(T)$
which is closed by hypothesis. It is easy to see that $\gamma (T)
= \gamma (B)$ and $\| B\| =  \|T\| $. Also, $B^\dag T$ is a
coisometry, with the same nullspace as $T$. So, by Lemma
\ref{L:coisometria}, $\gamma (B^\dag T P_\eme ) =\sang{N(T)}{
\eme}$. Now, using Lemma \ref{L:inversible} for $B $ and $B^\dag
TP_\eme $ and the fact that $BB^\dag TP_\eme = P_{R(T)}TP_\eme  =
TP_\eme $, we get
\[\gamma (T)\ \sang{N(T)}{ \eme} \leq \gamma (TP_\eme )\leq \|T\|
\ \sang{N(T)}{ \eme} ,
\]
because $R(B) =R( B^\dag)$, so that
$R( B^\dag T P_\eme ) \subseteq  R(B) =N(B)\orto $.
\end{proof}

\begin{rem}\rm
With the same ideas,  the following formulae
generalizing Lemma \ref{L:coisometria} and Proposition
\ref{P:gammas}, can be proved. 
\ben
\item
Let $U,V \in \ca$ be partial isometries. Then,
$\gamma(UV)=\sang{N(U)}{R(V)}$.
\item
If $A,B \in \ca $ have closed ranges, then
$$%\begin{equation}\label{desigualdad}
 \gamma(A)\gamma (B)\;\sang{N(A)}{R(B)} \leq \gamma (AB)\leq \|A\|\,\|B\|
 \;\sang{N(A)}{R(B)}  .
$$%\end{equation}
Note that the first inequality implies Proposition \ref{producto con rango cerrado}.

\een In particular, this gives the following formula for the sine
of an angle: given $\eme$ and $\ene$ two closed subspaces of
$\hil$, it holds $$ \sang{\ene}{ \eme} = \gamma
(P_{\ene\orto}P_\eme). $$ \EOE
\end{rem}

\subsection{Frames}

We introduce some basic facts about frames in Hilbert spaces. For
complete descriptions of frame theory and applications, the reader is referred to the
survey by Heil and Walnut \cite{[HeiWal]} or the books by Young \cite{[Y]} and Christensen
\cite{[liChr]}.
\begin{fed} \rm
Let $\H$ be a separable Hilbert space, and $\F = \fram$ a sequence in $\H$.
\ben
\item
$\F$ is called a \textsl{frame} if
there exist numbers $ A,B>0$ such that, for every $f \in \H$,
\beq\label{frame}
A\|f\|^2\leq \sum_{n\in \mathbb{N}}
|\pint{f , f_n}|^2 \leq B\|f\|^2
\end{equation}
\item The optimal constants  $A, B$ for equation (\ref{frame})  are called
the \textsl{frame bounds} for $\F$.
\item The frame $\F$ is called $tight$ if $A=B$, and $Parseval$ if
$A=B=1$.
\item Associated with $\F$ there exist an operator $T:\ell^2\to\hil$ such that 
$T(e_n)=f_n$ where $\bonn$ denotes the canonical basis of $\ell^2$. 
This operator is called the $synthesis$ operator of $\F$. For finite 
 frames we assume that the domain of the synthesis operator is 
$\mathbb{C}^m$ where $m$ is the number of vectors of the frame.
\een
\end{fed}

\begin{rem} \label{cosas} \rm
Let $\F= \fram$ be a frame in $\hil$ and $T$ its synthesis operator.
\ben
\item The frame bounds of  $\F$ can be computed in terms of the synthesis operator
\beq\label{bounds}
A=\gamma(T)^{2}  \quad \hbox { and } \quad B=\|T\|^2 .
\end{equation}

\item The adjoint $T^*\in  L(\hil , \ell^2 )$ of 
$T$, is   given by $\displaystyle
T^*(x  ) = \sum_{n\in\N} \api x  , f_n \cpi e_n$,
$x  \in \H$. It is called the \textsl{analysis operator} for $\F$.

\item The operator $S=TT^*$ is usually called frame operator and it is easy to see that 
\begin{equation}\label{formula pa el ese}
S f=\sum_{n\in \N} \pint{f , f_n}f_n \hspace{1cm} f\in \H.
\end{equation}
 It follows from (\ref{frame}) that  $A.I \leq S \leq B.I$, so that $S \in Gl( \hil )^+$. 
 Moreover, the optimal
constants  $A, B$ for equation (\ref{frame})  are $$ B = \|S\| =
\rho(S) \quad \hbox{ and } \quad A = \gamma (S) = \|S\inv \|\inv =
\min \{ \la : \la \in \sigma (S)\}. $$ 
Finally, from \eqref{formula pa el ese} we get
\[
f=\sum_{n\in \N} \pint{f , S\inv f_n}f_n \hspace{1cm} \forall\ f\in\hil.
\]
\item
The numbers $\{ \pint{f , S\inv f_n}\}$ are called the \textsl{frame
coefficients} of $f$. They have the following optimal property:
if $f = \sum_{n\in \N} c_n f_n$, for a
sequence $(c_n)_{n\in \N}$, then
$$
\sum_{n\in \N} |\pint{f , S\inv f_n}|^2 \le \sum_{n\in \N} |c_n|^2.
$$
The frame $\{S^{-1}f_n\}_{n\in\N}$ is called canonical dual frame. We shall return to dual frames in section 4.
\EOE
\een
\end{rem}

\section{Riesz frames and conditional Riesz frames.}

It was remarked by Christensen \cite {[liChr]} , p. 65, that given a frame
$\F = \fram$,  in practice it can be difficult to use the frame
decomposition $f = \sum \pint{f , S\inv f_n} f_n$ because it requires the 
calculation of  $S\inv$ or, at least, the frame coefficients 
$\pint{f , S\inv f_n}$. In order to get some of the advantages of
Riesz bases, Christensen introduced in \cite{[Chri1]} the
\textbf{projection method}, approximating $S$ and $S\inv$ by
finite rank operators, acting on certain finite dimensional spaces
$\H_n$ approaching $\H$. Later on, Christensen \cite {[Chri3]}
introduced two special classes of frames, namely \textbf{Riesz
frames} and \textbf{conditional Riesz frames}, which are well
adapted to some of these problems (see also \cite{[Caz1]}, \cite{[CazChri]}, and \cite{[CazChri2]}).

\medskip
\noi We need to fix some notations: Let  $\cB = \bonn$  be
%{\bf Notations}: Given
the canonical orthonormal basis of $\ell^2$ and  $I\subseteq \mathbb{N}$. \ben
\item
We denote $\ese_I = \gencl{e_n : n\in I}$ and $P_I = P_{\ese_I}$, the orthogonal projection onto
$\ese_I$.
\item
If  $I = \mathbb{I}_n := \{1,2,  \dots , n\}$,
we put $\ese_n$ for $\ese_I$.
\item
Given $\ene$ a closed subspace of  $\ell^2$,
we  denote  $\ene_n =\ene\cap \ese _n$,
$n \in \N$.
\item
If $\F = \fram $ is a frame for $\H$, we denote by
$\F_I = \{f_n\}_{n\in I}$.
\item
We say that $\F_I$ is a
\textsl{frame sequence} if it is a frame for $\gencl{\F_I}$. 
\item
$\F_I$ is called a \textsl{subframe} of $\F$ if it is itself 
a frame for $\H$.
\een

\bigskip

\noi Recall the definitions of Riesz frames and conditional Riesz frames.

\begin{fed} \rm
A frame $\F = \fram$ is called a \textsl{Riesz frame} 
if there exists $A, B > 0$ such that, for  every 
$I\subset \mathbb{N}$, the subfamily $\F_I$  is a frame
sequence with bounds $A, B$ (not necessarily optimal).

The sequence $\F$ is called a \textsl{conditional Riesz frame} 
if there are common bounds for the frame sequences $\F_{I_n}$, where
$\{I_n\}_{n=1}^\infty$ is a sequence of finite subsets of $\mathbb{N}$
such that $I_n \inc I_{n+1 }$ for every $n\in \N$ and $\displaystyle
\bigcup _{n\in \N} I_n = \N$.
\end{fed}

\begin{rem}\rm
Let $\F$ be a frame, and $T$ its synthesis operator. Given $I \inc \N$, then $\F_I$ is a frame sequence \sii
$R(TP_I ) $ is closed, and $\F_I $ is a subframe \sii $R(TP_I) =
\H$. In both cases the frame bounds for $\F_I$ are $A = \gamma
(TP_I )^2$ and $B = \|TP_I\|^2$. Using these facts
we get an equivalent definition of Riesz
frames: 
$\F$ is a Riesz frame if there exists $\eps>0 $ such that
$\gamma (TP_I)\ge \eps $ for every $I \inc \N$. \EOE
\end{rem}

\noi Proposition \ref{P:gammas} can be used to characterize 
Riesz frames in terms of the angles between the nullspace of the 
synthesis operator $T$ and the closed subspaces of
$\ell^2$ which are spanned by subsets of $\cB$.

\begin{pro}\label{P: equivalencia entre ker y Rframe}
Let $\F =\fram$  be a frame, and $T$ be its synthesis operator. Let $\ene=N( T)$. Then
$\F $ is a Riesz frame \sii
    \begin{equation}\label{condicion del ker en Rframe}
  c \ =\  \sup _{I \inc \N} \ \angf{\ene}{\ese_I}<1 . 
    \end{equation}
\end{pro}

\begin{proof}
By Proposition \ref{producto con rango cerrado},  $TP_I$ has closed range
iff $\angf{\ene}{\ese _I}<1$. By Proposition \ref{P:gammas},
$\gamma (TP_I)$ has an uniform
    lower bound \sii there exists a constant $c<1$ such that, for every
    $I \subseteq \mathbb{N}$, $\angf{\ene}{\ese_I}\leq c$.
\end{proof}

\begin{rem}\label{nucleo compatible}\rm
Let $\ene $ be a  closed subspace of $\ell^2$ and $\cB = \bonn$ be the canonical orthonormal
basis of $\ell^2$. If equation (\ref{condicion del ker en Rframe}) holds,  following the
terminology of \cite{[ACRS]}, we say that $\ene$ is ${\mathcal
B}$-$compatible$.\EOE
%One of the properties that a $\cB$- compatible subspace
%of $\hil$ must satisfy is the density of
%$\bigcup_{n=1}^{\infty}\ese_n$ in $\ese$ (see \cite{[ACRS]}).
%Actually, this result can be reformulated
%for a general family of subspaces, which contains  all %nullspaces of
%conditional Riesz frames.
\end{rem}

\bigskip

\noi
In the following Proposition, 
we state a characterization of $\mathcal{B}$- compatible subspaces of $\hil$,
proved in \cite{[ACRS]}.

\begin{pro}\label{P:equivalencia de compatibilidad de S}
Let $\ene$ be a closed subspace of $\ell^2$ and let
$\mathcal{B}=\bon$ be the canonical orthonormal basis of $\ell^2$.  For
$n\ \in \N$,  denote by  $\displaystyle
c_n = \sup_{J\inc \, \mathbb {I}_n} \angf{\ene_n}{\ese_J}$.
Then  the following conditions are equivalent:

\begin{enumerate}
    \item
    $\ene$ is $\mathcal{B}$-compatible.
    \item
    $\displaystyle c = \sup_{n\in \N}\ \angf{\ene}{\ese_n}<1$,
    and \ $ \displaystyle\sup_{n\in \N}\  c_{n} < 1$.
    \item
    \rm $\cl{\bigcup_{n\in \N}\ene_n}=\ene $ \it
    and \
   $ \displaystyle\sup_{n\in \N} \ c_n < 1$.
   \item
   There exists a constant $c<1$ such that $\angf{\ene}{\hil _I}\leq c$ for every
  finite subset   $I$ of $\N $ with $\ene\cap \ese_I=\{0\}$.
\QED
\end{enumerate}
\end{pro}

\noi
Proposition \ref{P:equivalencia de compatibilidad de S} can be
``translated" to frame language to get a characterization of Riesz frames, similar
to the one obtained by Christensen and Lindner in \cite{[ChrisLin]}:

\begin{teo}\label{CL}%[see Theorem 2.1 in \cite{[ChrisLin]}]
Let $\F = \fram$ be a frame and $T$ its synthesis operator.
Denote $\ene = N(T)$. Then the following conditions are equivalent:

\ben
  \item $\F$ is a Riesz frame.
    \item
    $\ene$ is $\mathcal{B}$-compatible.
  \item There exists an uniform lower
frame bound for every finite linearly independent frame sequence $\F_J$, $J\subset \mathbb{N}$.
\item There exists $d>0$ such that $\gamma (TP_J ) \ge d$, for every $J \in \N$ finite such that $\ene
\cap \ese_J = \{0\}$. \een
\end{teo}

\begin{proof}
If $I$ is a finite subset of $\N$ then $\ese_I\cap\ene=\{0\}$ if
and only if $\F_I$ is linearly independent. Then, conditions 3 and
4 are equivalent. By Propositions \ref{P:gammas} and \ref{P:equivalencia de
compatibilidad de S}, they are also equivalent
to the $\cB$-compatibility of $\ene$.

Suppose that
there exists a constant $d$ such that $0< d\leq \gamma (TP_{\ese_I})$ for every
finite subset $I \inc \mathbb{N}$ such that $\ese_I \cap \ene=\{0\}$. This
is equivalent to saying that there is a constant $c<1$ such that
$\angf{\ene}{\ese_I}\leq c$
for such kind of sets $I$. Using  Propositions \ref{P: equivalencia entre ker y
Rframe} and \ref{P:equivalencia de compatibilidad de S}, we conclude  that $\F$ 
is a Riesz frame. The converse is
clear.
\end{proof}

\noi
Now, we consider conditional Riesz frames. First of all, we state a result 
for this class of frames which is similar to Proposition \ref{P: equivalencia entre 
ker y Rframe}, and  whose proof  follows essentially the same lines.

\begin{pro}\label{P:equivcondicionalRF}
Let $\F = \fram$ and $\ene$ the nullspace of its synthesis operator. Then
 $\F$ is a conditional Riesz frame \sii there exists  a  sequence $\{I_n\}$ of
finite subsets of  \ $\mathbb{N}$
such that $I_n \inc I_{n+1 }$, %for every $n\in \N$,
 \begin{equation}
 \bigcup _{n\in \N} I_n = \N \quad \hbox{ and } \quad
  c \ = \ \sup_{n \in \N } \ \angf{\ene}{\ese_{I_n}}<1  \ , \ \ n\in \N .
    \end{equation}
\end{pro}

\bigskip

\noi
As a corollary of this Proposition we get the following result:

\begin{pro}\label{p:densidad de Sn en S}
Let $\F$  be a conditional Riesz frame, and $T$ its synthesis operator for
$\F$. Denote $\ene=N(T)$.
Then \rm
\begin{equation}\label{union}
\cl{\bigcup_{n=1}^{\infty}\ene_n}=\ene .
\end{equation}
\end{pro}

\noi
In order to prove this Proposition, we need the following technical Lemma.

\begin{lem}\label{union de los esen es densa en ese}
 Let $\ene$ be a closed subspace of $\ell^2$, a constant $c<1$ and a  sequence
 $\{I_n\}$ of finite subsets of $\mathbb{N}$ such that $I_n \inc I_{n+1 }$,
 $\bigcup _{n\in \N} I_n = \N$ and
$  \angf{\ene}{\ese_{I_n}}\leq c $, for every $n\in \N$.
Then \rm
\[
\cl { \ \bigcup_{n\in \N}\ene\cap \ese_{I_n} }=\ene \ .
\]
\end{lem}
\begin{proof}
Denote $Q_n = P_{I_n}$, $n \in \N$.
The assertion of the Lemma is equivalent to
$$
P_{\ene}\wedge Q_n \convsoti P_\ene  \ .
$$
Let $x \in\ell^2$ be a unit vector and let $\eps>0$. Let $k \in
\mathbb{N}$ such that $\displaystyle c^{2k-1}\leq\frac{\eps}{2}$.
By Proposition \ref{aproximacion a la interseccion},  for every
$n\geq 1$ it holds that
\[
\left\|\left(P_{\ene}Q_n\right)^k-P_{\ene}\wedge Q_n\right\|\leq
\frac{\eps}{2}.
\]
On the other hand, since $Q_nP_\ene\convsot P_\ene$ and the
function $f(x)=x^k$ is SOT-continuous on bounded sets (see, for example,
2.3.2 of \cite{[Pe]}),
there exists $n_0\geq 1$ such that, for every $n\geq n_0$,
\[
\left\|\left[\left(Q_n P_\ene\right)^k-P_{\ene}\right]\; x
\right\|<\frac{\eps}{2} \ .
\]
Then, for every $n\geq n_0$,
\begin{align*}
\left\|\left(P_\ene-P_{\ene}\wedge Q_n\right)\; x  \right\|&\leq
\left\|\left[P_\ene-\left(P_{\ene}Q_n\right)^k\right]\; x
\right\|+ \left\|\left(\left(P_{\ene}Q_n\right)^k-P_{\ene}\wedge
Q_n\right)\; x  \right\|<\eps  \ .
\end{align*}
\end{proof}

\medskip

\noi \textit{Proof of Proposition \ref{p:densidad de Sn en S}}.$\ $
Since $\F$ is a conditional Riesz frame, there exist $c<1$ and a  sequence $\{I_n\}$ of
finite subsets of  \ $\mathbb{N}$ such that $I_n \inc I_{n+1 }$,
$\bigcup _{n\in \N} I_n = \N$ and
$\angf{\ene}{\ese_{I_n}}\leq c $, for every $n\in \N$. By Lemma
\ref{union de los esen es densa en ese}, $\bigcup_{n\in \N}\ene\cap \ese_{I_n}$
is dense in $\ene$. Finally, for every $n \in \N$, there exists $m \in
\N$ such that $I_n \inc \mathbb {I}_m = \{ 1, 2,\dots , m\}$. Thus,
$\bigcup_{n\in \N}\ene\cap \ese_{I_n} \inc \bigcup_{m\in \N}\ene_{m}$.
$\ $\QED 

\bigskip

\noi As a consequence of Proposition \ref{p:densidad de Sn en S} we obtain the following Corollaries.

\begin{cor}
Let $\F $ be  a conditional Riesz frame with synthesis operator $T$ 
and suppose that $\dim N(T)<\infty$. Then $\F$ is a Riesz frame. 
Moreover, there exists $m \in \N$ such that $N(T) \inc \ese_m$.
\end{cor}

\begin{proof}
Denote by  $\ene = N(T)$. By Proposition \ref{p:densidad de Sn en
S}, $\ene$ satisfies equation (\ref{union}). Since $\dim
\ene<\infty$, then there exists $m \in \N$ such that $\ene=N(T)
\inc \ese_m$. Thus, in the terminology of Proposition
\ref{P:equivalencia de compatibilidad de S}, if 
$\displaystyle c_n = \sup_{J\inc \mathbb{I}_n} 
\angf{\ene_n}{\ese_J}$, then $c_n = c_m$ 
for every $n\ge m$. Therefore, by Proposition
\ref{P:equivalencia de compatibilidad de S}, $\F$ 
is a Riesz frame.
\end{proof}

\begin{cor}
Let $\F= \fram$ be a conditional Riesz frame. Given $n\in
\mathbb{N}$, denote by $S_n$ the frame operator of
$\{f_k\}^n_{k=1}$ and let $A_n$ be the minimum of the lower frame
bounds  of all frame subsequences of $\{S_n\mrai f_k\}_{k=1}^n$. If \
$\inf_n A_n > 0$, then $\F$ is a Riesz frame.
\end{cor}
\begin{proof}
Let $T$ be the synthesis operator of $\F$ and $\ene = N(T)$.
For each $n\in \mathbb{N}$, denote $\F_n = \{f_k\}^n_{k=1}$,
$\cB_n = \{ e_1, \dots , e_n\}$ and $P_n = P_{\ese_n}$. Note that
$TP_n:\ese_n\to \gencl{f_k : k=1,\ldots,n}$ can be considered, 
modulo an unitary operator, as the synthesis operator of $\F_n$. In this way, it holds that 
$S_n = TP_nT^*$. Also note that $\{S_n\mrai f_k\}_{k=1}^n$ is a Parseval 
frame, and $N(TP_n ) = N(S_n\mrai TP_n ) =
\ese \cap \hil_n = \ese_n$. So, by Lemma \ref{L:coisometria},
if $J\subset \{1,\ldots ,n\}$, the
lower frame bound $A_J$ of $\{S_n\mrai f_k\}_{k\in J}$ satisfies
$A_J=1-\angf{\ese_n}{\hil _J}^2$.
Using Propositions \ref{p:densidad de Sn en S} and 
\ref{P:equivalencia de compatibilidad de S},
the corollary follows.
\end{proof}

\subsection*{A counterexample}

The nullspace $\ene$ of the synthesis operator
of a conditional Riesz frame has the property of
``density": $\cl{\bigcup_{n=1}^{\infty}\ene_n}=\ene $,  where
$\ene_n$ is $\ene \cap \ese_n$. In the following example we show
that the converse is not true, i.e., we  construct a frame which is
not a conditional Riesz frame such that  its synthesis nullspace
$\ene$ satisfies $\cl{\bigcup_{n=1}^{\infty}\ene_n}=\ene $.

\medskip
\noi
We shall prove the assertion in an indirect way, by using Proposition
\ref{P:equivcondicionalRF} and the 
following fact:
if $\ene$ is a closed subspace of $\ell^2$ such that $\dim \ene
^\bot=\infty$, then there exists a frame $\F$ with 
synthesis operator $T$  such that $\ene =N(T)$.

\begin{exa}\label{Cassa}
Given $r>1$, if $\cB = \bonn$ denotes the canonical basis of $\ell^2$, let us
define the following orthogonal system:
\begin{align*}
x_1&=e_1-re_2+\frac{1}{r}e_3+\frac{1}{r^2}e_4+
\frac{1}{r^3}e_5+\frac{1}{r^4}e_6\\
x_2&=e_5-re_6+\frac{1}{r^5}e_7+\frac{1}{r^6}e_8+
\frac{1}{r^7}e_9+\frac{1}{r^8}e_{10}\\
&\vdots\\ x _n&=e_{4n-3}-re_{4n-2}+\frac{1}{r^{4n-3}}e_{4n-1}+
\frac{1}{r^{4n-2}}e_{4n}+\frac{1}{r^{4n-1}}e_{4n+1}+
\frac{1}{r^{4n}}e_{4n+2}.
\end{align*}
Let $\ene$ be the closed subspace generated by 
$\{ x _n \}_{n\in \mathbb{N}}$. By construction,
$\cl{\bigcup_{n=1}^{\infty}\ene_n}=\ene $. Moreover
$\{e_{4n-1}-re_{4n}: n\in \mathbb{N}\} \subset \ene ^\bot$, 
so $\dim \ene ^\bot =\infty$. 
By the remarks above,  there exists a frame $\F$ such that the nullspace of 
its synthesis operator is $\ene$.
We claim that this frame is not a conditional Riesz
frame. By Proposition \ref{P:equivcondicionalRF}, 
it suffices to verify that for every
sequence $J_1\subseteq J_2 \subseteq J_3\subseteq \ldots 
\subseteq J_n \nearrow \mathbb{N}$, it holds that
$\angf{\ene}{\ese_{J_k}}\convk 1$. Hence, fix such a 
sequence $\{J_k\}_{k \in \N}$ and take $0<\eps <1$.

\noi
Since $\displaystyle\|x _n \|^2 \le 1 +r^2 +  \frac{4}{ r^{8n-6}}$
for every $n \in \N$, there exists $n_0 \in \mathbb{N}$ such
that 
$$
1-\eps<\frac{1+r^2}{\|x _n \|^{2}} \hspace{1cm}\forall\ n\geq n_0.
$$ 
Note that, for $y \in \ene $ and $i \in \N$, if $\cM_i = \gen{e_{4i-3},e_{4i-2}}$, 
then 
\beq\label{mi}
\api y , x_i \cpi = 0 \iff  P_{\cM_i}y =0,
\end{equation}
because $P_{\cM_i}x_j \neq 0$ \iiff $j = i$.
Let $k\in \N$ be such that 
$$ 
j=\max \Big\{i\in \N :\;P_{\cM_i}
(\ene\cap\ese_{J_k})\neq 0 \Big\}\ge n_0. 
$$ 
By equation (\ref{mi}), $x_{h} \in (\ene \cap \ese_{J_k})^\bot$ for every $h>j$.
In particular, $x _{j+1} \in \ene\ominus (\ene \cap \ese_{J_k})$ and 
\[
1-\eps < \frac{1+r^2}{\|x _{j+1} \|^2} \leq
\frac{\|P_{J_k}x _{j+1}\|^2}{\| x _{j+1}\|^2}\leq \pint{
\frac{x _{j+1}}{\|x _{j+1}\|} \; , \; \frac{P_{J_k}x _{j+1}}{\|
P_{J_k}x _{j+1}\| }} \leq \angf{\ene}{\ese_{J_k}}
\]
A similar argument shows that
$1-\eps\leq\angf{\ene}{\ese_{J_m}}$, for every $m\geq k$. 
This implies that $\displaystyle
\liminf_{n\to\infty}\angf{\ene}{\ese_{J_n}}\geq 1-\eps$. Finally,
as $\eps$ is arbitrary, we get
$\angf{\ene}{\ese_{J_k}}\convk 1$. \EOE
\end{exa}

%\newpage

\section{Weighted dual frames.}

Let $\F= \fram$ be a fixed frame for a closed subspace 
$\cW $ of $\hil$ and let $\ese\subseteq \hil$ be another closed 
subspace such that $\hil=\cW \sudi   \ese^\bot$. As we have  mentioned  
in the introduction, an oblique dual frame of $\F$ in $\ese$ 
is a frame $\G= \framg$ for $\ese$ such that for every $f\in\cW$ it holds that

\begin{equation}\label{de nuevo la formula de los frames duales}
f=\sum_{n=1}^\infty\pint{f,\ g_n}f_n \hspace{1cm} \forall \ f\in\cW.
\end{equation}

\noi Such a dual frame has the minimal norm property if for 
every $f\in\cW$ the coefficients $\left\{\pint{f,\ g_n}\right\}_{n\in\N}$ 
have minimal $\ell^2$ norm. Christensen and Eldar proved in \cite{[ChriEl2]}
that the duals frames with the minimal norm property have the form
\begin{equation}\label{duales segun Christensen}
\framg=\left\{B(T^*B)^\dagger e_n\right\}_{n\in\N} \ , 
\end{equation}
where $\bonn$ denote the canonical orthonormal basis of $\ell^2$, and $B$ 
is a bounded operator with $R(B)=\ese$.

\noi On the other hand, let $\deledos$ be the set of all $D\in { Gl}(\ell^2)^+$ 
which are diagonal in the canonical basis $\bonn$. Each $D\in\deledos$  
defines an inner product $\Dpint{\cdot,\ \cdot}$ by means of
\[
\Dpint{x,\ y}=\pint{Dx,\ y} \ , \quad x, y \in \ell^2 \ .
\]
This inner product induces a \textit{weighted} norm $\|\cdot\|_D$ which is equivalent to the original one. 

\bigskip

\noi In this section, we are interested in dual frames such 
that their coefficients in the reconstruction formula 
\eqref{de nuevo la formula de los frames duales} minimize different weighted norms. 
We shall give explicit formulae for this class of dual frames that 
we call weighted dual frames. We also consider the particular 
case of weighted dual frames associated to a Riesz frame.

\bigskip

\noi First of all, let us recall some preliminary facts on generalized inverse s.

\begin{fed}\label{definicion de pseudoinversa} \rm
Given two Hilbert spaces $\hil$ and $\cK$, let $A\in L(\hil,\cK)$ be an 
operator with closed range. We say that $B\in L(\cK,\hil)$ is a generalized inverse  of $A$ 
if  $ABA=A$ and  $BAB=B$. 
\end{fed}

\begin{rems} Let $A\in L(\hil,\cK)$ with closed range, 
and let  $B\in L(\cK,\hil)$ be a generalized inverse  of $A$. Then 
\begin{enumerate}
	\item Both $AB$ and $BA$ are oblique projections, i.e. idempotent operators.       \item $R(B)$ is also closed.
	\item   The idempotent $AB$ and $BA$ induce decompositions of the Hilbert spaces $\hil$ and 
	$\cK$:  $\hil=N(A)\sudi R(B)$ and $\cK=R(A)\sudi N(B)$. 
	\item If $(AB)^*=AB$ and $(BA)^*=BA$, then $B$  is called the Moore-Penrose generalized inverse for $A$. 
	It is usually denoted by $A^\dagger$. In this case,  
	$AA^\dagger$ is the orthogonal projection onto $R(A)$ and 
	$A^\dagger A$ is the orthogonal projection onto $N(A)^\bot$. \EOE
	\end{enumerate}
\end{rems}

\noi Among the generalized inverses of an operator $A\in L(\ell^2,\hil)$, the following ones will be particularly important for us. In order to clarify the next statement, given a subspace $\ete$ of $\ell^2$ and $D\in\deledos$, the orthogonal complement of $\ete$ with respect to the the inner product $\Dpint{\cdot,\ \cdot}$ will be denoted by $\ete^{\bot_D}$.

\begin{lem}\label{gueited pseudoinversas}
Let $A\in L(\ell^2,\hil)$ be an operator with closed range, and $D\in\deledos$. Then, the operator $\Chi_D(A)=D^{-1/2}(AD^{-1/2})^\dagger$ is a generalized inverse  of $A$ such that
$\Chi_D(A)A$ is the orthogonal projection with respect to the weighted inner 
product $\Dpint{\, \cdot \, ,\,  \cdot \, }$ onto $N(A)^{\bot_D}$.
\end{lem}
\begin{proof}
Since $R(AD^{1/2})=R(A)$ it follows  that 
$$
A\ \Chi_D(A)\ A= P_{R(AD^{1/2})}A=A.
$$
On the other hand, 
\[
\Chi_D(A)\ A\ \Chi_D(A)=D^{-1/2}(AD^{-1/2})^\dagger A D^{-1/2}(AD^{-1/2})^\dagger=D^{-1/2}(AD^{-1/2})^\dagger=\Chi_D(A).
\]
Finally, some easy computation shows that an oblique projection $Q$ is $D$-orthogonal if and only if $DQ$ is selfadjoint. In our case
\begin{align*}
D\left(\Chi_D(A)A\right)&= D^{1/2}(AD^{-1/2})^\dagger A= D^{1/2}\Big(D^{-1/2}A^*(ADA^*)^\dagger\Big)A = A^*(ADA^*)^\dagger A,
\end{align*}
which is clearly selfadjoint. Therefore, $\Chi_D(A)A$ is a $D$-orthogonal projection and clearly $N\big{(}\Chi_D(A)A\big{)}=N(A)$.
\end{proof}

\noi Now, we are ready to give the explicit form of weighted dual frames.

\begin{pro}\label{gueited dual frames}
Let $\F= \fram$ be a fixed frame for a closed subspace $\cW$ of $ \hil$, $T$ 
its synthesis operator and let $\ese$  be another closed subspace of $\hil$ 
such that $\hil=\cW\sudi \ese^\bot$. 
Then, given $D\in\deledos$, the oblique dual frames such that for 
every $f\in\cW$ their coefficient in the reconstruction formula minimize the 
weighted norm $\|\cdot\|_D$ have the form
\[
\G = \framg=\{B(D^{1/2}T^*B)^\dagger D^{1/2} e_n\}_{n\in\N}
\]
where $\bonn$ denotes the canonical orthonormal basis of $\ell^2$ and $B\in L(\ell^2,\hil)$ is any operator with $R(B)=\ese$.
\end{pro}
\begin{proof}
Fix $B\in L(\ell^2,\hil)$ with range $\ese$ and let 
$\widehat{T}=B(D^{1/2}T^*B)^\dagger D^{1/2}$. First of all, note that 
$N(D^{1/2}T^*B)=N(B)$. So, $R(\widehat{T})=R(B)=\ese$ and therefore $\G$ is a frame. 

In order to prove that $\G$ is an oblique dual frame it is enough to prove that $T\widehat{T}^*$ is an oblique projection onto $\cW$. Actually, $T\widehat{T}^*$ is the projection onto $\cW$ parallel to $\ese^\bot$. Indeed, on one hand
\begin{align*}
(T\widehat{T}^*)^2 &= \Big(TD^{-1/2}(B^*TD^{-1/2})^\dagger B^*\Big)^2
= T D^{-1/2}\Big((B^*TD^{-1/2})^\dagger (B^*T D^{-1/2}) (B^*TD^{-1/2})^\dagger\Big) B^*\\
&=T \Big(D^{-1/2}(B^*TD^{-1/2})^\dagger B^*\Big)=(T\widehat{T}^*),
\end{align*}
which shows that $T\widehat{T}^*$ is a projection. On the other hand, since $N(D^{-1/2}(B^*TD^{-1/2})^\dagger)=R(B^*T)^\bot=R(B)^\bot$ and $R(D^{-1/2}(B^*TD^{-1/2})^\dagger)=R(T^*B)=R(T^*)$, it holds that $T\widehat{T}^*$ is the projection onto $\cW$ with nullspace $\ese^\bot$. 

Finally, in order to prove that the reconstruction coefficients minimize the 
weighted norm $\|\cdot\|_D$ we have to prove that 
$R(\widehat{T}^*)\subseteq N(T)^{\bot_D}$. 
%, where $N(T)^{\bot_D}$ denotes the orthogonal complement of $N(T)$ with respect to the inner product $\Dpint{\cdot,\ $\cdot}$. 
But, using the notation of Lemma \ref{gueited pseudoinversas}, we get   
$\widehat{T}^*T=\Chi_D(B^*T)B^*T$ and, therefore, using the same Lemma, 
$R(\widehat{T}^*)\subseteq N(B^*T)^{\bot_D}=N(T)^{\bot_D}$. $\ $
\end{proof}

\noi As we have already mentioned in the previous section, 
$\fram$ is a Riesz frames if and only if $N(T)$ is compatible with 
the canonical base (see Remark \ref{nucleo compatible}). 
If  $\pas{D}{\ene}$ denote the (unique) orthogonal projection onto the closed 
subspace $\ene$ of $\ell^2$ with respect to 
the inner product $\Dpint{\cdot,\ \cdot}$, it was 
proved in \cite{[ACRS]} that $\ene$ is compatible \sii 
\[
\sup_{D\in\deledos} \|\pas{D}{\ene}\|<\infty.
\]

\noi As a consequence of this result we obtain the following.

\begin{teo}
Let $\F= \fram$ be a frame for a closed subspace $\cW$ of $\hil$, $T$ 
its synthesis operator, $\ese$ another closed subspace of $ \hil$ 
such that $\hil=\cW\sudi \ese^\bot$ and $\G= \framg$ a fixed (sampling) frame 
for $\ese$ with synthesis operator $B$. Then, the following conditions are equivalent:
\ben
\item $\F$ is a Riesz frame on $\cW$. 
\item The oblique dual frames of 
$T$ with respect to $B$ that minimize the different weighted norms 
are bounded from above. In other words
\[
\sup_{D\in\deledos}\|B(D^{-1/2}T^*B)^\dagger D^{-1/2}\|<\infty \ .
\]
\een
\end{teo}
\begin{proof}
Fix $D \in \deledos$. 
We have already proved in Lemma \ref{gueited pseudoinversas} that 
\[
\Big(B(D^{-1/2}T^*B)^\dagger D^{-1/2}\Big)^*T = T^*B(D^{-1/2}T^*B)^\dagger 
D^{-1/2}=T^*B\ \Chi_D(T^*B)=1-\pas{D}{N(T)} \ .
\]
Hence 
\begin{align*}
\|B(D^{-1/2}T^*B)^\dagger D^{-1/2}\|&\leq \|B\|\ \|(D^{-1/2}T^*B)^\dagger D^{-1/2}\|\\
&=\|B\|\ \|(T^*B)^\dagger(T^*B)(D^{-1/2}T^*B)^\dagger D^{-1/2}\|\\
&\leq \|B\|\ \|(T^*B)^\dagger\|\ \|(T^*B)(D^{-1/2}T^*B)^\dagger D^{-1/2}\| \\
&=  \|B\|\ \|(T^*B)^\dagger\|\  \|1-\pas{D}{N(T)}\| \ ,
\end{align*}
and 
$$
 \|1-\pas{D}{N(T)}\| = \|T^*B(D^{-1/2}T^*B)^\dagger 
D^{-1/2}\| \le \|T ^* \|\  \|B(D^{-1/2}T^*B)^\dagger  D^{-1/2}\| \ .
 $$
 Therefore 
 $$
 \sup _{D \in \deledos}  \|B(D^{-1/2}T^*B)^\dagger  D^{-1/2}\| < \infty 
 \iff \sup _{D \in \deledos} \|1-\pas{D}{N(T)}\| < \infty \ , 
 $$
 which proves the Proposition.
\end{proof}

\end{document}